\documentclass[preprint,review,12pt]{elsarticle}

\usepackage{amssymb}

\journal{Journal}
\begin{document}

\begin{frontmatter}
\title{On solutions of some of unsolved problems in number theory, specifically on the distribution of primes}
\author{Ahmad Sabihi\fnref{fn1}} 
\ead{sabihi2000@yahoo.com}
\address{Teaching professor and researcher at some universities of Iran}
\fntext[fn1]{ Fourth Floor, Building 30, Poormehrabi Alley, Mofatteh St., Bozorgmehr Ave., Isfahan, Iran.} 
\begin{abstract}
We solve some famous conjectures on the distribution of primes. These conjectures are to be listed as Legendre's, Andrica's, Oppermann's, Brocard's, Cram\'{e}r's, Shanks', and five Smarandache's conjectures. We make use of both Firoozbakht's conjecture (which recently proved by the author) and Kourbatov's theorem on the distribution of and gaps between consecutive primes. These latter conjecture and theorem play an essential role in our methods for proving these famous conjectures. In order to prove Shanks' conjecture, we make use of Panaitopol's asymptotic formula for $\pi(x)$ as well.  
\end{abstract}

\end{frontmatter}

\section{Introduction}

Recently, the author proved Firoozbakht's conjecture\cite{Fe1},\cite{Fe2}.This conjecture plays an important role in proving most of the conjctures on the distribution of primes. In this paper, we show that this conjcture along with Kourbatov's theorem 1\cite{K}are really useful and powerful for our purpose. In Section 2, we prove Legendre's conjecture.  Legendre's conjecture states that there exists at least a prime number between $n^{2}$  and  $(n+1)^{2}$ for all natural numbers.If this conjecture is correct, the gap between any prime $p$ and the next largest prime would always be at most on the order $\sqrt{p}$ or gaps are $O(\sqrt{p})$.This conjecture has been recognized to have not been solved since over 200 years ago. In Section 3, we prove Andrica's conjecture \cite{An} in the two ways.Andrica's conjecture states that the inequality $\sqrt {p_{n+1}}-\sqrt {p_{n}}<1$ holds for all $n$, where $p_{n}$ is the $n$th prime number.In Section 4, Oppermann's conjecture  is proven for every $n>1$. This conjecture is one of the unsolved problems in number theory, specifically on the distrbution of primes and was proposed by mathematician Ludvig Oppermann in 1882\cite{OP}.Oppermann's conjecture states that there is at least one prime as $p_{1}$ and one prime as $p_{2}$ so that 
\begin{equation}
n^{2}-n<p_{1}<n^{2}<p_{2}<n^{2}+n
\end{equation} 
for natural numbers $n\geq 2$. If the conjecture is true, then the largest possible gaps between two consecutive prime numbers could be at most proportional to twice the square root of numbers.In Section 5, Brocard's conjecture using the proven Oppermann's conjecture is proved.The conjecture says us that there exist at least four primes between $(p_{n})^{2}$  and  $(p_{n+1})^{2}$ for $n>1$, where $p_{n}$ is the $n$th prime number. In Section 6, we make a proof for Cram\'{e}r's conjecture.Cram\'{e}r's conjecture\cite{C}states that gaps between consecutive prime numbers can have a supermum 1 with regard to ${(\log p_{n})^{2}}$ ($\log$ refers to natural logarithm throughout the paper) as ${\lim_{n \rightarrow \infty}sup\frac{(p_{n+1}-p_{n})}{(\log p_{n})^{2}}=1}$.In Section 7, an easily proof of Shanks' conjecture is made.Shanks' conjecture \cite{S} (( ${p_{n+1}-p_{n})\sim (\log p_{n})^{2}}$) gives a somewhat stronger statement than Cram\'{e}r's. In Sections 8 to 12, we make the proofs of the first, second,third and fifth Smarandache's conjctures\cite{Sm1},\cite{Sm2}using proven Firoozbakht's conjcture,Kourbatov's theorem 1, and proven Andrica's conjecture and a disproof of his fourth conjecture in some special cases. These conjctures generalize Andrica's conjcture and will be  discussed in detail in their related Sections.  
\setcounter{equation}{0}
\section{Legendre's conjecture} \label{sc2}

As we should know, Legendre's conjecture states that there exists at least a prime number between  $n^{2}$  and  $(n+1)^{2}$ for all natural numbers.   

\textbf{\textit{Proof}}

According to the proven Firoozbakht's conjecture \cite{Fe1} and Kourbatov's theorem1 \cite{K}
\begin{equation}
{p_{k+1}-p_{k}<(\log p_{k})^{2}-\log p_{k}-1}      
\end{equation}
 ${for~~k>9~~~or~~~p_{k}\geq p_{10}=29}$
 
Thus,
\begin{equation}
{p_{k+1}-p_{k}<(\log p_{k})^{2}-\log p_{k}-1<(\log p_{k})^{2}}      
\end{equation}
 ${for~~k>9~~~or~~~p_{k}\geq p_{10}=29}$

Let $p_{k}$ be the greatest prime number right before $n^{2}$, then  $p_{k+1}$ should be between $n^{2}$  and  $(n+1)^{2}$.

Contradiction, assume there is no such $p_{k+1}$ between them,then  $p_{k}<n^{2}$ and $p_{k+1}>(n+1)^{2}$ 

In such a case, 

$\log p_{k}<2\log n$ and so $(\log p_{k})^{2}<4(\log n)^{2}$  and $p_{k+1}-p_{k}>2n+1$   

This means that 
\begin{equation}
{2n+1<p_{k+1}-p_{k}<(\log p_{k})^{2}<4(\log n)^{2}}      
\end{equation}
${for~~k\geq 11~~~or~~~n\geq 6}$

Trivially, the inequality (2.3) does not hold since $2n+1>4(\log n)^{2}$ for $n\geq 11$ and this implies that $p_{k+1}$ must be between $n^{2}$  and  $(n+1)^{2}$ and Legendre's conjecture would be true for all $n\geq 11$. On the other hand, this conjecture is also correct for $n\leq 10$, hence it holds for all $n\geq 1$. 
\setcounter{equation}{0}
\section{Andrica's conjecture}\label{sc3}
Andrica's conjecture states that the inequality $\sqrt {p_{n+1}}-\sqrt {p_{n}}<1$ holds for all $n$, where $p_{n}$ is the $n$th prime number. If we manipulate the inequality, it changes to 
\begin{equation}
{p_{k+1}-p_{k}<2\sqrt {p_{k}}+1}      
\end{equation}
\textbf{\textit{Proof}}

As is mentioned in Section2, regarding (2.2) and (3.1), we should prove
\begin{equation}
{p_{k+1}-p_{k}<(\log p_{k})^{2}<2\sqrt {p_{k}}+1}      
\end{equation}
\textbf{\textit{The first solution}}

Let $p_{k}$ be replaced by $x \in \mathbb R$, then we show $(\log x)^{2}<2\sqrt {x}+1$ for $x\geq 121$.
Let $y=2\sqrt {x}-(\log x)^{2}+1$ be a function of variable $x$ defined for $x\geq 121$.
$y(121)=0.000393$. Easily, we prove the derivation of $y$ is positive for all $x\geq 121$, i.e $y'>0$. 
\begin{equation}
{y'=\frac{\sqrt {x}-2\log x}{x}}      
\end{equation}
and $y'(121)=0.0116$.
Just, we show that the numerator (3.3) i.e ${\sqrt {x}-2\log x}>0$ for $x\geq 121$. Again, let $z=\sqrt {x}-2\log x$ and $z'=\frac{\sqrt {x}-4}{2x}$. Therefore, $z'>0$ for $x>16$, and so $z'(x)>0$ for $x\geq 121$ and $z(x)>z(121)>0$ for all $x\geq 121$,then $y'>0$ and $y>0$ for all $x\geq 121$ and the inequality $(\log p_{k})^{2}<2\sqrt {p_{k}}+1$ certainly   holds for $p_{k}\geq 121$.

Andrica's conjecture also holds for all $p_{k}<121$. Hence, it holds for all $k$.

\textbf{\textit{The second solution}}

We show that, if we replace $p_{k}$ by positive integer, $n$, in (3.2)
\begin{equation} 
(\log n)^{2}<2\sqrt {n}+1
\end{equation}
for $n\geq 190$

Easily,$(\log n)^{2}-1<2\sqrt {n}$, then $(1-\frac{1}{(\log n)^{2}})<\frac{2\sqrt{n}}{(\log n)^2}$

Taking $(\log n)^{2}$th power,
\begin{equation} 
(1-\frac{1}{(\log n)^{2}})^{(\log n)^{2}}<(\frac{2\sqrt{n}}{(\log n)^2})^{(\log n)^{2}}
\end{equation}
Trivially, analogous with $(1-\frac{1}{n})^{n}<\frac{1}{e}$  for $n\geq 1$, we have $(1-\frac{1}{(\log n)^{2}})^{(\log n)^{2}}<\frac{1}{e}$  for $n>e$.

Easily, we show $(\frac{2\sqrt{n}}{(\log n)^2})^{(\log n)^{2}}>\frac{1}{e}$ or $(\frac{\sqrt{n}}{(\log n)^2})>\frac{1}{2}e^{-\frac{1}{(\log n)^{2}}}$ for $n\geq 190$.
Since $\frac{1}{2}e^{-\frac{1}{(\log n)^{2}}}<\frac{1}{2}$ for all natural numbers and $\frac{\sqrt{190}}{(\log 190)^2}=0.50066>0.5$ we should prove that $\frac{\sqrt{n}}{(\log n)^2}>\frac{\sqrt{190}}{(\log 190)^2}$ for $n>190$.

This means that we should prove that the sequence $\frac{\sqrt{n}}{(\log n)^2}$ is strictly increasing for $n\geq 190$. A simple calculation shows that the sequence is increasing for all natural numbers 190 to 320. We only need to show it is correct for $n\geq 321$. We show that the inequality
\begin{equation} 
\frac{\sqrt{n+1}}{(\log (n+1))^2}>\frac{\sqrt{n}}{(\log n)^2}
\end{equation}
holds for $n\geq 321$.

Manipulating (3.6),
$\sqrt{n+1}~(\log n)^2>\sqrt{n}~(\log (n+1))^2$ and $\sqrt{1+\frac{1}{n}}>(1+\frac{\log (1+\frac{1}{n})}{\log n})^{2}$

Taking $n$th power,
\begin{equation}
(1+\frac{1}{n})^{\frac{n}{2}}>(1+\frac{\log (1+\frac{1}{n})}{\log n})^{2n}
\end{equation}
then
\begin{equation}
\{(1+\frac{1}{n})^{n}\}^{\frac{1}{2}}>\{(1+\frac{\log (1+\frac{1}{n})}{\log n})^{\frac{\log n}{\log (1+\frac{1}{n})}}\}^{\frac{2n\log (1+\frac{1}{n})}{\log n}}
\end{equation}
Trivially, $\{(1+\frac{1}{n})^{n}\}^{\frac{1}{2}}>2^{\frac{1}{2}}$ for $n\geq 1$. Thus, we need to prove that 
\begin{equation}
\{(1+\frac{\log (1+\frac{1}{n})}{\log n})^{\frac{\log n}{\log (1+\frac{1}{n})}}\}^{\frac{2n\log (1+\frac{1}{n})}{\log n}}<2^{\frac{1}{2}}
\end{equation}
for $n\geq 321$.
Trivially, $\{(1+\frac{\log (1+\frac{1}{n})}{\log n})^{\frac{\log n}{\log (1+\frac{1}{n})}}\}^{\frac{2n\log (1+\frac{1}{n})}{\log n}}<e^{\frac{2}{\log n}}<2^{\frac{1}{2}}$ for $n\geq 321$.

Therefore, (3.9), (3.8),(3.7), and consequently (3.6) hold for $n\geq 321$ and (3.4) holds for $n\geq 190$ or $p_{k}\geq 190$. 

Andrica's conjecture also holds for all $p_{k}<190$. Hence, it holds for all $k\geq 1$.
\setcounter{equation}{0}
\section{Oppermann's conjecture}\label{sc4}
Oppermann's conjecture states that there is at least one prime as $p_{1}$ and one prime as $p_{2}$ so that 
\begin{equation}
n^{2}-n<p_{1}<n^{2}<p_{2}<n^{2}+n
\end{equation} 
for natural numbers $n\geq 2$.

\textbf{\textit{Proof}} 

Regarding (2.2), Let $p_{k}$ be the greatest prime right before $n^{2}-n$, then $p_{k+1}$ should be between $n^{2}-n$ and $n^{2}$ .

Thus,
\begin{equation}
p_{k}<n^{2}-n
\end{equation}
Assume  $p_{k+1}$ does not exist between $n^{2}-n$ and $n^{2}$, then
\begin{equation}
p_{k+1}>n^{2}
\end{equation}
From (4.2) and (4.3), 

$p_{k+1}-p_{k}>n$   and 
\begin{equation}
n<p_{k+1}-p_{k}<(\log p_{k})^{2} 
\end{equation}
On the other hand, from (4.2)

$\log p_{k}<\log(n^{2}-n)$ and so 
\begin{equation}
(\log p_{k})^{2}<(\log n+\log(n-1))^{2}
\end{equation}
${for~~k>9~~~or~~~p_{k}\geq p_{10}=29}$. If $29\leq p_{k}<n^{2}-n$, then $n\geq 6$.  

Trivially,  $\log n+\log(n-1)<2\log n$ and $(\log n+\log(n-1))^{2}<4(\log n)^{2}$

Just, we prove that 
\begin{equation}
(\log n+\log(n-1))^{2}<4(\log n)^{2}<n 
\end{equation}
for $n\geq 75$.
 
Consider $2\log n<\sqrt{n}$ for $n\geq 75$  

Let $y=\sqrt{x}-2\log x$, then $y'=\frac{\sqrt{x}-4}{2x}$ which implies $y'>0$ for $x>16$. Also, we know that $y>0$ for $x\geq 75$. Thus, (4.6) holds for $n\geq 75$.

Therefore, holding (4.4),(4.5), and (4.6) leads us to a contradiction and our assumption,which asserts $p_{k+1}$ does not exist between $n^{2}-n$ and $n^{2}$ is incorrect for $n\geq 75$. This means that Oppermann's conjecture is true for all  $n\geq 75$. Oppermann's conjecture  trivially holds for $2\leq n<75$ and consequently holds for $n\geq 2$.

The second part of Oppermann's conjecture also holds easily and similarly with suppostion: Let $p_{k}$ be the greatest prime right before $n^{2}$, then $p_{k+1}$ should be between $n^{2}$ and $n^{2}+n$. 

Assume  $p_{k+1}$ does not exist between $n^{2}$ and $n^{2}+n$, then $p_{k+1}>n^{2}+n$.

Similarly, we have
\begin{equation}
n<{p_{k+1}-p_{k}<(\log p_{k})^{2}-\log p_{k}-1}<(\log p_{k})^{2}<4(\log n)^{2}      
\end{equation}
 ${for~~k>9~~~or~~~p_{k}\geq p_{10}=29}$
 
where leads us to $n<4(\log n)^{2}$. But,this is a contradiction since $n>4(\log n)^{2}$ for $n\geq 75$. This means that $p_{k+1}$ exists between $n^{2}$ and $n^{2}+n$ and Oppermann's conjecture holds for $n\geq 75$ and consequently for $n\geq 2$.    
\setcounter{equation}{0}
\section{Brocard's conjecture}\label{sc5}
The conjecture says us that there exist at least four primes between $(p_{n})^{2}$  and  $(p_{n+1})^{2}$ for $n>1$, where $p_{n}$ is the $n$th prime number.
 
\textbf{\textit{Proof}}
 
The proof is easily made by proven Oppermann's conjecture in Section 4.
We decompose the gap between $(p_{n})^{2}$ and $(p_{n+1})^{2}$ into the four segments,
\begin{itemize}
\item The gap between $(p_{n})^{2}$ and $p_{n}(p_{n}+1)$
\item The gap between $p_{n}(p_{n}+1)$ and $(p_{n}+1)^{2}$ 
\item The gap between $(p_{n}+1)^{2}$ and $(p_{n}+1)(p_{n}+2)$
\item The gap between $(p_{n}+1)(p_{n}+2)$ and $(p_{n}+2)^{2}$ 
\end{itemize}
We only need to prove that 
\begin{equation}
(p_{n+1})^{2}\geq (p_{n}+2)^{2}
\end{equation}
Let $(p_{n+1})^{2}-(p_{n})^{2}=(p_{n+1}-p_{n})(p_{n+1}+p_{n})$. Trivially, $(p_{n+1}-p_{n})\geq 2$ for $n\geq 2$, thus $p_{n+1}+p_{n}\geq 2p_{n}+2$ so

$(p_{n+1})^{2}-(p_{n})^{2}=(p_{n+1}-p_{n})(p_{n+1}+p_{n})\geq 4p_{n}+4$. Hence

$(p_{n+1})^{2}\geq (p_{n})^{2}+4p_{n}+4=(p_{n}+2)^{2}$

Therefore, there exists at least a prime number between each of the above four gaps and Oppermann's conjecture is proved for $n\geq 2$.
\setcounter{equation}{0}
\section{Cram\'{e}r's conjecture}\label{sc6}
This conjecture states
\begin{equation}
{\lim_{n \rightarrow \infty}sup\frac{(p_{n+1}-p_{n})}{(\log p_{n})^{2}}=1}
\end{equation}

\textbf{\textit{Proof}}

As mentioned in the previous conjectures, regarding (2.2) we have

\begin{equation}
\frac{p_{n+1}-p_{n}}{(\log p_{n})^{2}}<1
\end{equation}
for $n>9$.

This immediately implies (6.1).

Also, we have
\begin{equation}
p_{n+1}-p_{n}=O((\log p_{n})^{2})
\end{equation}
The inequality (6.3) shows us that for sufficiently large integers,$n$, we should have 
\begin{equation}
|{p_{n+1}-p_{n}}| \leq A|(\log p_{n})^{2}|
\end{equation}
We easily find that (2.2) implies Cram\'{e}r's conjecture with $A=1$ and hence 
\begin{equation}
p_{n+1}-p_{n}<(\log p_{n})^{2}
\end{equation}
for when $n$ tends to infinity and Cram\'{e}r's conjecture is satisfied. 
\setcounter{equation}{0}
\section{Shanks' conjecture}\label{sc7}
Shanks improved Cramer's conjecture by another strongly conjecture. He states that
\begin{equation}
({p_{n+1}-p_{n})\sim (\log p_{n})^{2}}
\end{equation}
for sufficiently large integers or when $p_{n}$ tends to infinity.

\textbf{\textit{Proof}}

The proof is easily made by proven Firoozbakht's conjecture \cite{Fe1}, Kourbatov's theorem 1 \cite{K} (the relation (2.1)) and Panaitopol's theorem 1 \cite{P}. Panaitopol's theorem 1 states that
\begin{equation}
\pi(x)=\frac{x}{\log x-1-\frac{k_{1}}{\log x}-\frac{k_{2}}{(\log x)^{2}}-...-\frac{k_{n}(1+\alpha_{n}(x))}{(\log x)^{n}}} 
\end{equation}
where $k_{1},k_{2},...,k_{n}$ are given by the recurrence relation
$k_{n}+1!k_{n-1}+2!k_{n-2}+...+(n-1)!k_{1}=n.n!,~~n=1,2,3...$
and $\lim_{x \rightarrow \infty}\alpha_{n}(x)=0$ or $\alpha_{n}(x)=O(\frac{1}{\log x})$.
We easily check using $\pi(x)$ given by (7.2) and letting $x=p_{k}$ that
\begin{equation}
k<\frac{\log p_{k}}{\log p_{k+1}-\log p_{k}}<\frac{ p_{k}}{\log  p_{k}-1-\frac{k_{1}}{\log  p_{k}}-\frac{k_{2}}{(\log  p_{k})^{2}}-...-\frac{k_{n}(1+\alpha_{n}( p_{k}))}{(\log  p_{k})^{n}}-|O(\sqrt{\log p_{k}})|} 
\end{equation}
The following inequality is known to be true
\begin{equation}
\log(x+y)-\log x<\frac{y}{x}~~~for~~every~~ x,y>0
\end{equation}
Let $y=p_{k+1}-p_{k}$ and $x=p_{k}$ into the relation (7.4) and combine to (7.3)
\begin{eqnarray}
\frac{(\log  p_{k})^{2}-\log  p_{k}-k_{1}-\frac{k_{2}}{\log  p_{k}}-...-\frac{k_{n}(1+\alpha_{n}( p_{k}))}{(\log  p_{k})^{n-1}}-|O(\sqrt{\log p_{k}})|\log p_{k}}{p_{k}}<\log p_{k+1}-\log p_{k}\nonumber\\<\frac{p_{k+1}-p_{k}}{p_{k}}~~~~~~~~~~~~~~~~~~~~~~~~~~~~~~~~~~~~~~~~~~~~~~~~~~~~~
\end{eqnarray}
and gives us
\begin{equation}
p_{k+1}-p_{k}>(\log  p_{k})^{2}-\log  p_{k}-k_{1}-\frac{k_{2}}{\log  p_{k}}-...-\frac{k_{n}(1+\alpha_{n}( p_{k}))}{(\log  p_{k})^{n-1}}-|O(\sqrt{\log p_{k}})|\log p_{k}
\end{equation}
Combining (7.6) to Kourbatov's theorem 1 gives us
\begin{eqnarray}
(\log  p_{k})^{2}-\log  p_{k}-k_{1}-\frac{k_{2}}{\log  p_{k}}-...-\frac{k_{n}(1+\alpha_{n}( p_{k}))}{(\log  p_{k})^{n-1}}-|O(\sqrt{\log p_{k}})|\log p_{k}<p_{k+1}-p_{k}\nonumber\\<(\log p_{k})^{2}-\log p_{k}-1~~~~~~~~~~~~~~~~~~~~~~~~~~~~~~~~~~~~~~~~~~~\end{eqnarray}
${for~~k>9~~~or~~~p_{k}\geq p_{10}=29}$

Dividing both sides by $(\log  p_{k})^{2}$ and tending $p_{k}$ to infinity, we have
\begin{equation}
1<\lim_{p_{k} \rightarrow \infty} \frac{p_{k+1}-p_{k}}{(\log  p_{k})^{2}}<1
\end{equation}
This means that
\begin{equation} 
\lim_{p_{k} \rightarrow \infty} \frac{p_{k+1}-p_{k}}{(\log  p_{k})^{2}}=1 
\end{equation}
and Shanks' conjecture is proven.
\setcounter{equation}{0}
\section{First Smarandache's conjecture \cite{Sm2}}\label{sc8}
This conjecture says us that equation $(p_{n+1})^{x}-(p_{n})^{x}=1$, where $p_{n}$ and $p_{n+1}$ denote the $n$th and $(n+1)$th primes respectively, has a unique solution for $0.5<x\leq 1$.
The maximum solution occurs for $n=1$, i.e. $3^{x}-2^{x}=1$ when $x=1$. The minimum solution occurs for $n=31$, i.e. $127^{x}-113^{x}=1$ when $x=0.567148...=a_{0}$.

\textbf{\textit{Proof}}

The proof is comprised of the three steps.
\begin{enumerate}
\item There is no solution for $x>1$

Let $x=1+\epsilon$, where $\epsilon >0$,then

$(p_{n+1})^{1+\epsilon}-(p_{n})^{1+\epsilon}=(p_{n+1}-p_{n})\{(p_{n+1})^{\epsilon}+(p_{n+1})^{\epsilon -1}p_{n}+...+(p_{n})^{\epsilon}\}$.
Since, $p_{n+1}-p_{n}\geq 2$ for $n\geq 2$ and $\{(p_{n+1})^{\epsilon}+(p_{n+1})^{\epsilon -1}p_{n}+...+(p_{n})^{\epsilon}\}>1$ for $\epsilon >0$, implies  $(p_{n+1})^{1+\epsilon}-(p_{n})^{1+\epsilon}\geq 2$ for $n\geq 2$ and $\epsilon >0$. This means that we showed $(p_{n+1})^{x}-(p_{n})^{x}\neq 1$. 
\item There is no solution for $x\leq 0.5$ 

According to Andrica's theorem (Section 3), $(p_{n+1})^{\frac{1}{2}}-(p_{n})^{\frac{1}{2}}<1$. We should show that 
\begin{equation}
(p_{n+1})^{\frac{1}{t}}-(p_{n})^{\frac{1}{t}}<(p_{n+1})^{\frac{1}{2}}-(p_{n})^{\frac{1}{2}}<1
\end{equation} 
for $t>2$ , $t\in \mathbb R$.

Let $y=(p_{n+1})^{\frac{1}{t}}-(p_{n})^{\frac{1}{t}}-1$
and $z=(p_{n})^{\frac{1}{t}}$, then $z'=\frac{-(p_{n})^{\frac{1}{t}}\log p_{n}}{t^{2}}$.

$y'=\frac{(p_{n})^{\frac{1}{t}}\log p_{n}-(p_{n+1})^{\frac{1}{t}}\log (p_{n+1})}{t^{2}}<0$, since $(p_{n})^{\frac{1}{t}}\log p_{n}<(p_{n+1})^{\frac{1}{t}}\log (p_{n+1})$

for $t\geq 2$.

This means that $y<0$ is a strictly decreasing function for $t\geq 2$. This implies that  function $y$ has no any solution for $x=\frac{1}{t}\leq 0.5$.
\item We found out that $y=(p_{n+1})^{\frac{1}{t}}-(p_{n})^{\frac{1}{t}}-1$  is a continuously and strictly decreasing function for all $t>0$. As we showed $y>0$ for $0<t<1$ i.e. $\frac{1}{t}=x=1+\epsilon>1$, also $y<0$ for $t\geq 2$ or $x\leq 0.5$. We therefore find out due to having contiuously and strictly decreasing property of $y$ for all real numbers $t>0$, it must be zero $y=0$ for a unique value $x$ based upon the intermediate value theorem \cite{Si}.
\end{enumerate}
\setcounter{equation}{0}
\section{Second Smarandache's conjecture \cite{Sm2}}\label{sc9}
The conjecture generalizes Andrica's conjecture ($A_{n}=(p_{n+1})^{\frac{1}{2}}-(p_{n})^{\frac{1}{2}}<1$) to $B_{n}=(p_{n+1})^{a}-(p_{n})^{a}<1$ ,where $a<a_{0}$.

\textbf{\textit{Proof}}

 We should show that for $a<a_{0}$ 

$(p_{n+1})^{a}<(p_{n})^{a}+1$ , then taking $a$-th root, $p_{n+1}<((p_{n})^{a}+1)^{\frac{1}{a}}=p_{n}+\frac{1}{a}((p_{n})^{a})^{\frac{1}{a}-1}+...+1$

Therefore, we show 
\begin{equation}
p_{n+1}-p_{n}<\frac{1}{a}(p_{n})^{(1-a)}+...+1
\end{equation}
for $a<a_{0}$ 

Regarding (2.2), we would show that
\begin{equation}
p_{n+1}-p_{n}<(\log p_{n})^{2}-\log p_{n}-1<(\log p_{n})^{2}<\frac{1}{a}(p_{n})^{(1-a)}+...+1      
\end{equation}
${for~~n>9}$ and $a<a_{0}$

For our purpose, it is sufficient that we only prove
\begin{equation}
(\log p_{n})^{2}<\frac{1}{a}(p_{n})^{(1-a)}      
\end{equation}
$for~~n>9$ and $a<a_{0}$ since $\frac{1}{a}(p_{n})^{(1-a)}<\frac{1}{a}(p_{n})^{(1-a)}+...+1$.        

Let $p_{n}$ be replaced by $x\in \mathbb R$, then for real numbers $x\geq 5850$ and $a<a_{0}$, we should have   
\begin{equation}
(\log x)^{2}<\frac{1}{a}x^{(1-a)}      
\end{equation}
Let 
\begin{equation}
y=\frac{1}{a}x^{(1-a)}-(\log x)^{2}
\end{equation}
and
\begin{equation}
y'=\frac{(1-a)}{a}x^{-a}-\frac{2}{x}\log x
\end{equation}
Just we want to show that $y>0$ and $y'>0$ for $x\geq 5850$ and $a=a_{0}$. Certainly, if we have the result for when $a=a_{0}$,we will also have it for all $a<a_{0}$ since $\frac{1}{a}x^{(1-a)}>\frac{1}{a_{0}}x^{(1-a_{0})}$ and $\frac{(1-a)}{a}x^{-a}>\frac{(1-a_{0})}{a_{0}}x^{-a_{0}}$ for $a<a_{0}$.

For $a=a_{0}$, (9.5) and (9.6) are obtained
\begin{equation}
y=1.76320819x^{0.432852}-(\log x)^{2}
\end{equation}
and
\begin{equation}
y'=0.76320819x^{-0.567148}-\frac{2}{x}\log x
\end{equation}
Checking for $x=5850$ implies $y\geq 0.08077$ and $y'>0$

Manipulating the inequality $y'>0$ defining $y'$ in (9.8), we should show 
\begin{equation}
\frac{x}{(\log x)^{2.3095}}>9.33 
\end{equation}
for $x>5850$

Let $z(x)=\frac{x}{(\log x)^{2.3095}}-9.33$, then $z'=\frac{(\log x)^{1.3095}\{\log x-2.3095\}}{(\log x)^{4.619}}$. 

Easily checking gives us 

$z'>0$ for $x\geq 5850>e^{2.3095}$ and so $z(x)>z(5850)=30.5>0$ for $x\geq 5850$. Therefore, (9.9) holds and consequently $y'>0$ defining $y'$ in (9.8) and $y>0$ defining $y$ in (9.7) for $x\geq 5850$

This means that (9.4) holds for $a\leq a_{0}$ and
\begin{equation}
(\log p_{n})^{2}<\frac{1}{a_{0}}(p_{n})^{(1-a_{0})}<\frac{1}{a}(p_{n})^{(1-a)}<\frac{1}{a}(p_{n})^{(1-a)}+...+1
\end{equation}
holds for all $p_{n}>5850$ and the inequalities (9.2) and (9.1) hold for $p_{n}>5850$. Therefore, the conjecture holds for $p_{n}>5850$. Trivially by calculating, the conjecture holds for $p_{n}<5850$ and finally holds for $p_{n}\geq 2$. 
\setcounter{equation}{0}
\section{Third Smarandache's conjecture \cite{Sm2}}\label{sc10}
This conjecture generalizes Andrica's conjecture ($A_{n}=(p_{n+1})^{\frac{1}{2}}-(p_{n})^{\frac{1}{2}}<1$) to $C_{n}=(p_{n+1})^{\frac{1}{k}}-(p_{n})^{\frac{1}{k}}<\frac{2}{k}$ ,where $k\geq 2$.

\textbf{\textit{Proof}}

Arguing similarly to second Smarandach's conjecture
\begin{equation}
(p_{n+1})^{\frac{1}{k}}<(p_{n})^{\frac{1}{k}}+\frac{2}{k} 
\end{equation}
Taking $k$th power
\begin{equation}
p_{n+1}<((p_{n})^{\frac{1}{k}}+\frac{2}{k})^{k}=p_{n}+2(p_{n})^{(\frac{k-1}{k})}+...+(\frac{2}{k})^{k}
\end{equation}
for  $k\geq 2$.

Thus, we expect to have
\begin{equation}
p_{n+1}-p_{n}<2(p_{n})^{(\frac{k-1}{k})}+...+(\frac{2}{k})^{k}
\end{equation}
for  $k\geq 2$.

Regarding (2.2), it is sufficient to show
\begin{equation}
(\log p_{n})^{2}<2(p_{n})^{(\frac{k-1}{k})}+...+(\frac{2}{k})^{k}
\end{equation}
and regarding proven Andrica's conjecture in Section 3, we showed that
\begin{equation}
(\log p_{n})^{2}<2\sqrt {p_{n}}+1 
\end{equation}
for $p_{n}\geq 121$. 
 
Therefore, it is easily verifiable that
\begin{equation}
2\sqrt {p_{n}}+1<2(p_{n})^{(\frac{k-1}{k})}+2(\frac{k-1}{k})(p_{n})^{(\frac{k-2}{k})}...+(\frac{2}{k})^{k}
\end{equation}
for $k\geq 2$ and $p_{n}\geq 121$ since $\frac{k-1}{k}\geq \frac{1}{2}$ and $\frac{k-2}{k}\geq 0$. This means that (10.5),(10.4),(10.3), and (10.2) hold and consequently (10.1) holds for $k\geq 2$ and $p_{n}\geq 121$. Investigating this conjecture for $p_{n}<121$ shows that it is correct for all $n$ and for $k\geq 2$ 
\setcounter{equation}{0}
\section{Fourth Smarandache's conjecture \cite{Sm2}}\label{sc11}
This conjecture would also generalize Andrica's conjecture to 
\begin{equation}
D_{n}=(p_{n+1})^{a}-(p_{n})^{a}<\frac{1}{n}
\end{equation} 
where $a<a_{0}$ and $n$ big enough, $n=n(a)$, holds for infinitely many consecutive primes.
 
\textbf{\textit{Disproof}}

This conjecture cannot be correct for sufficiently large integers,$n$, with constant value $a$. This is because of if $n$ tends to infinity and $a=cte.$, then $(p_{n+1})^{a}-(p_{n})^{a}<\frac{1}{n}$ is not correct since 
\begin{equation}
{\lim_{n \rightarrow \infty}\{(p_{n+1})^{a}-(p_{n})^{a}}\}<{\lim_{n \rightarrow \infty}\frac{1}{n}=0 }
\end{equation}   
This means that 
\begin{equation}
{\lim_{n \rightarrow \infty}(p_{n+1})^{a}}<{\lim_{n \rightarrow \infty}(p_{n})^{a}}
\end{equation}
Taking $a$-th root gives us
\begin{equation}
{\lim_{n \rightarrow \infty}p_{n+1}}<{\lim_{n \rightarrow \infty}p_{n}}
\end{equation} 
which leads to a contradiction.

If we even tend $n$ to infinity and $a$ to zero simultaneously in inequality $(p_{n+1})^{a}-(p_{n})^{a}<\frac{1}{n}$ depending on how tends each of $(p_{n+1})^{a}-(p_{n})^{a}$  and $\frac{1}{n}$ to zero,the result may be correct or not. Therefore, one is not able to make decision on the result.

a) Is this still available for $a_{0}<a<1$?

According to the previous argument, it is not correct since $n$ tends to infinity.

b) Is there any rank $n_{0}$ depending on $a$ and $n$ such that (11.1) is verified for all $n\geq n_{0}$?

This may be correct if we take $a$ as a sufficiently small value.  
\setcounter{equation}{0}
\section{Fifth Smarandache's conjecture \cite{Sm2}}\label{sc12}
This conjecture says us that inequality $\frac{p_{n+1}}{p_{n}}\leq \frac{5}{3}$ holds for all $n$ and the maximum occurs at $n=2$.

\textbf{\textit{Proof}}

 Trivially, this conjecture is verified for $n=1,2,3$. The proven Firoozbakht's conjecture for all $n$ implies that 

\begin{equation}
\frac{p_{n+1}}{p_{n}}<(p_{n})^{\frac{1}{n}}
\end{equation} 

Considering the inequality (12.1) for $n\geq 4$, we verify it for $n=4$ and $\frac{p_{5}}{p_{4}}=1.571...<(p_{4})^{\frac{1}{4}}=1.6266..<\frac{5}{3}$

Easily, we check
\begin{equation}
\frac{p_{5}}{p_{4}}<(p_{4})^{\frac{1}{4}}<\frac{5}{3}
\end{equation} 
\begin{equation}
\frac{p_{6}}{p_{5}}<(p_{5})^{\frac{1}{5}}<(p_{4})^{\frac{1}{4}}<\frac{5}{3}
\end{equation} 
and finally conclude that
\begin{equation}
\frac{p_{n+1}}{p_{n}}<(p_{n})^{\frac{1}{n}}<(p_{n-1})^{\frac{1}{(n-1)}}<...<(p_{4})^{\frac{1}{4}}<\frac{5}{3}
\end{equation} 
which gives us
\begin{equation}
\frac{p_{n+1}}{p_{n}}<\frac{5}{3} 
\end{equation}   
for $n\geq 4$ and completes the proof.\\\\

\textbf{Acknowledgment}

The author would like to thank mathematicians Farideh Firoozbakht for introducing her strong conjcture and Alexi Kourbatov for presenting his nice and impressive theorem.I am indebted to them for making use of their brilliant thoughts.I also thank Prof. Smarandache for his nice generalizations on Andrica's conjecture 

\end{document}